\documentclass{article}
\def\var{\varepsilon}
\def\R{{{\rm I}\!{\rm R}}}

\def\H{\hbox{\bf H}}

\def\m{{\bf m}}

\usepackage{mathtools}
\usepackage{amssymb}
\usepackage{amsfonts}
\usepackage[us,12hr]{datetime}
\usepackage{fancyhdr}
\usepackage[active]{srcltx}

\usepackage{pdfsync} 
\DeclareSymbolFontAlphabet{\mathbb}{AMSb}
\usepackage{latexsym}
\usepackage{amsmath}
\usepackage{amsfonts}
\usepackage{mathbbol}
\usepackage{srcltx}
\usepackage {amsbsy}
\usepackage{amssymb}
\usepackage{epsfig}
\usepackage{latexsym}
\usepackage{amsmath}
\usepackage{amsfonts}

\newcommand{\eproof} {\framebox{}}

\newcounter{mycounter}[section]

\begin{document}

\bigskip\bigskip\bigskip
\begin{center} {\Large \bf
{ A magneto-viscoelasticity problem with a singular memory kernel}
 }
\end{center}
\normalsize \vspace{0.2cm}

\medskip
\begin{center}
{\bf Sandra Carillo$^{1,2}$, Michel Chipot$^{3}$, Vanda Valente$^{4}$  and Giorgio Vergara Caffarelli
}\footnote{\begin{enumerate}
\item Dipartimento {\it Scienze di Base e Applicate
    per l'Ingegneria},       \textsc{Sapienza} - Universit\`a di Roma,  16, Via A. Scarpa, 00161 Rome, Italy
\item ~I.N.F.N. - Sezione Roma1, Gr. IV - Mathematical Methods in NonLinear Physics,  Rome, Italy 
\item Universit\"at Zurich, Institute of Mathematics, Zurich, Switzerland
\item Istituto per le Applicazioni del Calcolo, C.N.R, via dei Taurini 19, 00185 Roma, Italy
\end{enumerate}
}
\end{center}

\parindent 0em

\medskip
\begin{abstract}
The existence  of solutions to a one-dimensional  problem arising in
magneto-viscoelasticity is here considered. Specifically, a non-linear system of integro-differential 
equations is analyzed; it is obtained coupling an integro-differential equation modeling the viscoelastic behaviour,  in which the kernel represents the relaxation function, with  the 
non-linear partial differential equations modeling the presence of a magnetic field. 
The case under investigation generalizes a previous study since the relaxation 
function is allowed to be unbounded at the origin, provided it belongs to $L^1$;
 the magnetic model equation adopted, as in the previous results 
 \cite{GVS2010, GVS2012, CSVVumi,CSVV} is  the penalized  Ginzburg-Landau  
 magnetic evolution equation.
\end{abstract}
\section{Introduction}
\setcounter{equation}{0}
The study of magneto-viscoelastic materials is motivated by the interest on mechanical properties of innovative 
materials widely studied in a variety of applications. In particular, as far as the coupling between mechanical 
and magnetic effects is concerned, the interest is motivated by  new materials such as Magneto Rheological Elastomers  or, 
in general, magneto-sensitive polymeric composites (see  
\cite{Hossain-et-al-2015a, Hossain-et-al-2015b}  and references therein). 
A variational approach to study multiscale models, in this context, is given in \cite{Ethiraj-et-al2015}. 
The results here presented are connected to a wide research project concerning the analytical study of 
differential and integro-differential models connected  to mechanical properties of materials. Thus, in 
\cite{CSVVumi, CSVV, VV} magneto-elasticitcity problems are considered,  in \cite{GVS2010, GVS2012} 
magneto-viscoelasticity problems are studied. Then, in turn, the case of a $1$-dimensional, and of a 
$3$-dimensional, body is investigated under the assumption of a regular kernel representing the relaxation 
modulus.  
Later,  materials with memory characterized by a singular kernel integro-differential  equations are studied  in 
\cite{[81], CVV2013a, S2014}. Indeed, as pointed out therein, the case is of interest not only to model different 
physical behaviours but also under the analytical viewpoint. The interest in {\it singular kernel} problems goes 
back to {Boltzmann} \cite{Boltzmann} and later, is testified, analytically, by the results of {Berti} \cite{Berti}, 
{Giorgi and Morro} \cite{Giorgi-Morro92}, Grasselli and Lorenzi \cite{Grasselli-Lorenzi91}  and {Hanyga} et al.
\cite{Hanyga, HanygaS, Hanyga-MS2}.  In addition,  {\it fractional derivative} models, since the works of  
{Rabotnov} \cite{Rabotnov} and Koeller \cite{Koeller}, are employed in  \cite{Adolfsson et al, Enelund et al, Enelund-Olsson}.
Here, a viscoelastic body is studied under the assumption of a  relaxation 
modulus, modeled by a $L^1$ function, coupled with a magnetic field.

The problem to study is concerned with the behaviour of a viscoelastic body subject 
also to  the presence of a magnetic field. The body is assumed to be one-dimensional. 
In particular the problem under investigation is motivated by a great interest in the realization 
of new materials which, on one side couple a viscoelastic behaviour with a magnetic one, see 
{\begin{equation} \label{eq-evo}
\left \{ 
\begin {array}{l} \displaystyle u_t({t})-\int_0^t {G}(t-\tau) u_{xx}(\tau) d\tau - u_1 
-\int_0^t \frac{\lambda}{2}(\Lambda(\m)\cdot \m)_{x} d\tau\ =\int_0^t  f(\tau) d\tau\, 
\\
\\
 \displaystyle\mathbf m_t+{\mathbf m}\frac{|\mathbf m|^2-1}
{\delta}+\lambda\Lambda (\mathbf m) u_{x}-\mathbf m_{xx}=0,
\end{array} \!\!\! \text{ in }   {\cal Q} \right.
\end{equation}}
together with the initial and boundary conditions 
\begin{equation} \label{ic}
u(\cdot,0)=u_0=0,\,\,\,
 \mathbf
m(\cdot,0)={\m}_0,\quad |{\m}_0|=1\quad {\rm in}\,\, \Omega\,,
\end{equation}
\begin{equation} \label{bc}
u=0,\qquad  \frac{\partial \mathbf m}{\partial \mbox{\boldmath
$\nu$}}=0 \qquad \mbox {on}\quad \Sigma=\partial\Omega \times
(0,T)\,,
\end{equation}
where $\Omega =(0,1)$, ${\cal Q}:=\Omega \times (0,T)$  and  
${{\mathbb{\cal M}}}\equiv(0,\m)$, letting $\m=(m_1,m_2)$, is the magnetization vector, 
orthogonal to the
conductor so that, since ${\bf u} \equiv (u,0,0)$,  
when both quantities are written in $\R^3$;  in addition, {\boldmath $\nu$} is the outer unit normal at the boundary 
$\partial \Omega$, $\Lambda$ is a
linear operator defined by $\Lambda(\m)=(m_2,m_1)$, the scalar
function $u$ is the displacement in the direction of the conductor itself, here identified with the $x-$axis and
$ \lambda$ is a positive parameter. In addition, the term $f$
represents an external force which also includes the deformation
history.\\
Moreover we assume:
 \begin{equation}
u_1 \in L^2(\Omega), \quad \m_0 \in \H^1(\Omega),~~ f\in L^2({\cal Q}).
\end{equation}

The model adopted here  to describe  the magneto-elastic interaction  is introduced in
 \cite{He}, \cite{CSVVumi}, \cite{CSVV},
\cite{VV} and the case of magneto-viscoelastic  regular behaviour is given in \cite{GVS2010, GVS2012}.

In fact, the kernel in the linear integro-differential equation, which represents the relaxation function $G$,  is assumed here  to satisfy weaker functional requirements with respect to the {\it classical} regularity requirements. 
   In particular, the relaxation function $ {G}(t)$ is assumed to 
   be such that 
    \begin{equation}\label{newG}
    G \in L^1(0,T)\cap C^2(0,T), 
    ~~ \forall T\in \R^+~;
    \end{equation}
    the relaxation function $G(t)$ is assumed to satisfy the further requirements,
   which follow from the physics of the model,    
   \begin{equation}\label{G}
    G(t)>0,\qquad \dot{G}(t)\le 0, \qquad \ddot{G}(t)\ge 0,\qquad
    t\in(0,\infty).
    \end{equation}
    Note that, in the  {\it classical} model the relaxation function, 
  further to satisfy    conditions (\ref{G}), 
     is assumed to be $  C^2[0,T],    ~~ \forall T\in \R^+$.
To this aim, in the following Section 2, a suitable  sequence of approximated classical 
problems
  is constructed. In the same Section also some apriori estimates  are obtained.
  Crucial in our analysis is the assuption $u_0=0$.

The  subsequent Section 3, is devoted to prove the  existence of a weak solution to the 
problem (\ref{eq-evo}) with the initital and boundary conditions (\ref{ic}) - (\ref{bc}). 

\section{Approximated problems and a priori  estimates}

In this Section, the approximation strategy is devised and, then, some estimates which are 
needed to prove the existence results are given.

First of all,  observe that, the reason why equation (\ref{eq-evo})$_1$ is written under the 
form of an  evolution equation is that the {\it classical} model \cite{D,D1} is not defined since it depends on 
$G(0)$ and on the integral of $\dot G$ which is not assumed to be in $L^1$ at the origin. 
However, in our case,
even if the the kernel $G$ of the integral equation is singular at the origin, 
 the  regularity requirements it is supposed to satisfy  
(Cf. (\ref{newG})), guarantee that the  {\it classical}  problem can be adopted 
to model  the   magneto-viscoelastic behaviour of the material as soon 
as we consider a time $t>0$. Hence,
 here a sequence of {\it time-translated} approximated problems is constructed.
 Specifically, let ${\varepsilon}$ denote  
    a {\it small parameter} $0<{\varepsilon}\ll1$ and consider an approximated 
    problem corresponding to each value of the parameter \,${\varepsilon}$ defined via a ${\varepsilon}$ time translation, that is, let us introduce, corresponding to each ${\varepsilon}>0$,
the {\it translated  relaxation function} $ G^{\varepsilon}(\cdot):=G({\varepsilon}+\cdot)$. 
Furthermore, it is coupled with a penalized version of the magnetization equation with  {\it penalization parameter} $0<{\delta}\ll1$, i.e. adopting  the same model in \cite{GVS2010}
where, now, the magnetization problem is coupled with a translated viscoelasticity equation. 
Then, we can introduce  the problem $P^{\varepsilon}$ given by : 
\begin{equation} \label{eqso}
\displaystyle{
\left \{ \begin {array}{l} \displaystyle u^{\varepsilon}_{tt}- G^{\varepsilon}(0)u^{\varepsilon}_{xx}-
\int_0^t \dot{G}^{\varepsilon}(t-\tau) u^{\varepsilon}_{xx}(\tau) d\tau
-\frac{\lambda}{2}(\Lambda(\m^{\varepsilon})\cdot \m^{\varepsilon})_{x}=f
\\
\\
 \displaystyle\mathbf m^{\varepsilon}_t+{\mathbf
m^{\varepsilon}}\frac{|\mathbf m^{\varepsilon}|^2-1}
{\delta}+\lambda\Lambda (\mathbf m^{\varepsilon}) u^{\varepsilon}_{x}-\mathbf m^{\varepsilon}_{xx}=0, 
\end{array}\right. ~\text{ in } {\cal Q}}
\end{equation}
together with the initial and boundary conditions
\begin{equation} \label{ic-eps-del}
u^{\varepsilon}(\cdot,0)=u_0=0,\,\,\, u^{\varepsilon}_t(\cdot,0)=u_1,\quad \mathbf
m^{\varepsilon}(\cdot,0)={\m}_0,\quad {\rm in}\,\, \Omega\,,
\end{equation}
\begin{equation} \label{bc-eps-del}
 u^{\varepsilon}=0,\qquad 
 \frac{{\partial \mathbf m}^{\varepsilon}}{\partial \mbox{\boldmath
$\nu$}}=0 \qquad \mbox {on}\quad \Sigma=\partial\Omega \times (0,T)\,,
\end{equation}

where $G^{\varepsilon}(t)\in  C^2[0,T]$ and, hence, the non-linear integro-differential problem $P^{\varepsilon}$ 
is well defined. Specifically, according to \cite{GVS2010}, such a problem admits a {\it unique strong solution}. \medskip

\noindent
{\bf Lemma 2.1}
 \label{2.1}  { \it
 Let $\bar{u}$ denote a solution to the problem 

\begin{equation}\label{eql2}
\displaystyle \bar{u}_{tt}- G^{\varepsilon}(0)\bar{u}_{xx}-
\int_0^t \dot{G}^{\varepsilon}(t-\tau) \bar{u}_{xx}(\tau) d\tau = F \text{ in } {\cal Q},
\end{equation} where in the r.h.s. $F\in L^2({\cal Q})$.
The initial and boundary conditions, in turn, are
\begin{equation} \label{ic+bc-G}
\bar u(\cdot,0)=u_0=0,\,\,\, \bar u_t(\cdot,0)=u_1,\qquad {\text{\rm in }} \Omega 
\end{equation}
\begin{equation} \label{ic+bc-G} \bar u=0,\quad \qquad\qquad  {\text{\rm on }}
\quad \Sigma=\partial\Omega \times (0,T)\,,
\end{equation}
  it follows
  {
\begin{equation}\label{ineq-lemma1} \displaystyle 
{1\over{2}} \int_{\Omega} G(t+\varepsilon)\, \vert\bar{u}_x\vert^2\, dx + {1\over{2}}
\int_{\Omega} \vert\bar{u}_t\vert^2\, dx
\le  {1\over{2}} \int_{\Omega} \vert{u}_1\vert^2\, dx ds  +{1\over{2}}
\int_{\Omega} G(\varepsilon)\, \vert{u}_x(0)\vert^2+
\end{equation}  }
$$+ \int_{\Omega} \int_0^t   F\,\bar{u}_t\, dx \,ds ~.$$}\rm

\medskip
\noindent
{\bf Proof.} Here, for the reader convenience, we give the proof which follows the original one by Dafermos \cite{D, D1}.
First of all, equation  (\ref{eql2}), when we change the integration variable $\tau$ into $s=t-\tau$ add and subtract the term 
  $$ \int_0^t \dot{G}^{\varepsilon}(s) u_{xx}(t) ds= \left[G^{\varepsilon}(t)-G^{\varepsilon}(0)\right]  u_{xx}(t)~,$$
can be written in the following equivalent form
{\large\begin{equation} \label{eql1}
\displaystyle{   \bar{u}_{tt} - {G}(t+\varepsilon)\bar{u}_{xx} + \int^t_0
  \dot{{G}}(s+\varepsilon)\left[\bar{u}_{xx}(t)-\bar{u}_{xx}(t-s)\right] ds = F.
  }  
\end{equation}  }
Then  multiplication of   equation (\ref{eql1}) by $\bar{u}_t$,
and integration over $\Omega$   gives 
{
\begin{equation}\label{form2}\begin{array}{cl@{\hspace{0.5ex}}c@{\hspace
{1.0ex}}l} \displaystyle{ {1\over{2}} {d\over{dt}}\int_{\Omega}\vert\bar{u}_t\vert^2 dx
\,+\,\int_{\Omega}{G}(t+\varepsilon) \bar{u}_x\, \bar{u}_{xt} dx \,+}
\\   \displaystyle{ \!\!\!\!\!\!\!\!\!\!\!\!\!\!\!\!
+\int_{\Omega} \bar{u}_t(t) \, dx\, \int_0^t \dot{{G}}(s+\varepsilon)
\,\left[\bar{u}_{xx}(t)-\bar{u}_{xx}(t-s)\right]ds  = \int_{\Omega}
F \bar{u}_t\, dx}
\end{array}\end{equation}  }
that is, since $\bar{u}_{t}$ is independent of $s$

{\begin{equation}\label{form3}\begin{array}{cl@{\hspace{0.5ex}}c@{\hspace
{1.0ex}}l}   

\displaystyle{{1\over{2}} {d\over{dt}}\int_{\Omega}\vert\bar{u}_t\vert^2 dx
\,+\,{1\over{2}} {d\over{dt}} \int_{\Omega}{G}(t+\varepsilon) \vert\bar{u}_x\vert^2\,dx
-{1\over{2}}  \int_{\Omega} \dot{G}(t+\varepsilon) \vert\bar{u}_x\vert^2\,dx \,-}
\\ \\ \displaystyle{\!\!\!\!\!\!\!\!
-\int_{\Omega}\, dx \, \int_0^t \dot{{G}}(s+\varepsilon) \bar{u}_{xt}
\,\left[\bar{u}_{x}(t)-\bar{u}_{x}(t-s)\right] ds = \int_{\Omega}
F \bar{u}_t\, dx .}
        \end{array}\end{equation}   }
Now, observe that {
\begin{eqnarray*}
\label{form3b}
 &
\displaystyle{-\int_{\Omega} \int_0^t \dot{G}(s+\varepsilon) \, \bar{u}_{xt}
 [\bar{u}_x(t)-\bar{u}_x(t-s)] \, dx ds }\\
&
\displaystyle{
= -{1\over{2}} {d\over{dt}}\int_0^t ds \int_{\Omega}
\dot{G}(s+\varepsilon)\vert\bar{u}_x(t)-\bar{u}_x(t-s)\vert^2  dx  }
\\
&
 \displaystyle{
+{1\over{2}} \int_{\Omega} \dot{G}(t+\varepsilon) \vert\bar{u}_x(t)-\bar{u}_x(0)\vert^2
 dx  }\\
&
\qquad\qquad\displaystyle{
-\int_{\Omega} \int_0^t \dot{G}(s+\varepsilon)
\bar{u}_{xt}(t-s)[\bar{u}_x(t)-\bar{u}_x(t-s)] \, dx ds } \\
&
\displaystyle{= -{1\over{2}} {d\over{dt}}\int_0^t ds \int_{\Omega}
\dot{G}(s+\varepsilon)\vert\bar{u}_x(t)-\bar{u}_x(t-s)\vert^2  dx }\\
&
\qquad\qquad\qquad\qquad\displaystyle{+ {1\over{2}}
\int_{\Omega} \dot{G}(t+\varepsilon) \vert\bar{u}_x(t)-\bar{u}_x(0)\vert^2
 dx   }
 \\ 
    &
    \displaystyle{\qquad\qquad+
    \int_{\Omega} \int_0^t \dot{G}(s+\varepsilon) \frac{d}{ds}\bar{u}_x (t-s)
    [\bar{u}_x(t)-\bar{u}_x(t-s)] \, dx ds }\\
    &\displaystyle{= -{1\over{2}} {d\over{dt}}\int_0^t ds \int_{\Omega}
    \dot{G}(s+\varepsilon)\vert\bar{u}_x(t)-\bar{u}_x(t-s)\vert^2 \, dx  }\\
    &\qquad\displaystyle{+ {1\over{2}}
    \int_{\Omega} \dot{G}(t+\varepsilon) \vert\bar{u}_x(t)-\bar{u}_x(0)\vert^2
    \, dx  }
     \\  
    &  \qquad\qquad\qquad\qquad  \displaystyle{-
    {1\over{2}}\int_{\Omega} \int_0^t \dot{G}(s+\varepsilon) \frac{d}{ds}
    \vert\bar{u}_x(t)-\bar{u}_x(t-s)\vert^2 \, dx ds }\\
\end{eqnarray*}

\begin{eqnarray*}
\label{form3b}
&\displaystyle{= -{1\over{2}} {d\over{dt}}\int_0^t ds \int_{\Omega}
\dot{G}(s+\varepsilon)\vert\bar{u}_x(t)-\bar{u}_x(t-s)\vert^2 \, dx  }\\
&\qquad\qquad\displaystyle{ +
{1\over{2}}\int_{\Omega} \int_0^t \ddot{G}(s+\varepsilon)
\vert\bar{u}_x(t)-\bar{u}_x(t-s)\vert^2 \, dx \,ds.}
\end{eqnarray*}}
Substitution of the latter  in (\ref{form3}),  
combined with integration  over the time interval $(0,t)$,  taking into
account the sign conditions (\ref{G}),  implies (\ref{ineq-lemma1}) {
and, hence, in the case of homogeneous  initial displacement condition, 
completes the proof.
\hfill$\Box$

\bigskip
This last estimate,   later on, is combined with the following one.

\medskip
\noindent
{\bf Lemma 2.2} 
 \label{2.2}  { \it
 Let $({u}^{\varepsilon}, \m^{\varepsilon})$ denote a solution to the problem 
{\rm (\ref{eqso})--(\ref{bc-eps-del})}, then it follows that we have  {
\begin{equation*}
\!\displaystyle{{1\over{2}}
\int_{\Omega} G^{\varepsilon}(t)\vert{u}^{\varepsilon}_x\vert^2\, dx + {1\over{2}}
\int_{\Omega} \vert{u}^{\varepsilon}_t\vert^2\, dx + \int_0^t \int_{\Omega}  \vert{\m}^{\varepsilon}_t\vert^2\, dx  + {1\over{2}}
\int_{\Omega} \vert{\m}^{\varepsilon}_x\vert^2\, dx + 
 }\end{equation*}
 \begin{equation}\label{ineq-lemma2a}
\displaystyle{ \frac{\lambda}{2} \int_{\Omega}
 \Lambda(\m^{\varepsilon})\cdot \m^{\varepsilon}
{u}^{\varepsilon}_x\, dx +{1\over{4}}
\int_{\Omega} {(\vert{\m}^{\varepsilon}\vert^2-1)^2\over \delta}\, dx\le
}
\end{equation}

{\begin{equation*}
\displaystyle{
\le  \int_0^t \int_{\Omega} f u^{\varepsilon}_t dx +
+ {1\over{2}} 
\int_{\Omega} {\vert{\m}_{0,x}\vert^2}\, dx +
 {1\over{2}} \int_{\Omega} \vert{u}_1\vert^2\, dx ~~.}
\end{equation*}}  }
}\rm

{\bf Proof} \quad
For sake of  simplicity, all the superscripts are omitted.
Consider the second equation in (\ref{eqso})
\begin{equation}
 \displaystyle\mathbf m_t+{\mathbf
m}\frac{|\mathbf m|^2-1}
{\delta}+\lambda\Lambda (\mathbf m) u_{x}-\mathbf m_{xx}=0.\
\end{equation}
Taking the scalar product with $\m_t$, after integration over $\Omega$, it
follows {
\begin{equation}\label{interm-m1}
\!\!\displaystyle{ \int_{\Omega}  \vert{\m}_t\vert^2\, dx  + {1\over{4}} {d\over{dt}}
\int_{\Omega}{\left( \vert{\m}\vert^2-1\right)^2 \over \delta}\, dx + {1\over{2}} {d\over{dt}}
\int_{\Omega} {\vert{\m}_x\vert^2}\, dx + \lambda \int_{\Omega} (\Lambda (\m)\cdot \m_t )u_x =0
}\end{equation}  }
and hence, after integration over $(0,t)$, recall \eqref{ic}, i.e. $ |{\m}_0|=1$,
\begin{equation}
\label{latter}
\!\!\displaystyle{ \int_0^t \int_{\Omega}  \vert{\m}_t\vert^2\, dx  + {1\over{4}} 
\int_{\Omega}{\left( \vert{\m}\vert^2-1\right)^2 \over \delta}\, dx + {1\over{2}} 
\int_{\Omega} {\vert{\m}_x\vert^2}\, dx + \lambda \int_0^t \int_{\Omega} (\Lambda (\m)\cdot \m_t) u_x dx}= 
\end{equation}
$$\displaystyle{ {1\over{2}} \int_{\Omega} {\vert{\m_0}_x\vert^2}\, dx}
$$}

Now, mulplying the first equation in (\ref{eqso}) by $u_t$ and integrating over $\Omega$, recalling Lemma 2.1,
it follows
 {
\begin{equation}\label{ineq1}\begin{array}{cl@{\hspace{0.5ex}}c@{\hspace
{1.0ex}}l}
\!\!\!\!\!\!\!\! \displaystyle{{1\over{2}}   \int_{\Omega} G(t+\varepsilon)\, \vert{u}_x\vert^2\, dx + 
{1\over{2}} \int_{\Omega} \vert{u}_t\vert^2\, dx}\le
\\ \\
\!\!\!\!\!\!\!\!\!\!\!\! \displaystyle{\le \int_{\Omega} \int_0^t \left[f + \frac{\lambda}{2}(\Lambda(\m)\cdot \m)_{x}\right]\,{u}_t\, dx \,ds + {1\over{2}}
\int_{\Omega} G(\varepsilon)\, \vert{u}_x(0)\vert^2\, dx + {1\over{2}} \int_{\Omega} 
\vert{u}_1\vert^2\, dx},
\end{array}\end{equation}
where, since  the initial homogeneous datum $(u_0=0)$ is assigned,  the term where
$G(\varepsilon)$ appears cancels.
            Since $$(\Lambda (\m)\cdot \m_t)={1\over{2}}{d\over{dt}}(\Lambda (\m)\cdot \m)~~ \text{and}~~(\Lambda (\m)\cdot \m_x)={1\over{2}}{d\over{dx}}(\Lambda (\m)\cdot \m),$$ 
the combination of  (\ref{latter}) and (\ref{ineq1})  allows to write
\begin{equation}\label{ineq-lemma2}\begin{array}{cl@{\hspace{0.1ex}}c@{\hspace
{0.1ex}}l}
\!\!\!\!\!\!\!\! \displaystyle{{1\over{2}}   \int_{\Omega} G(t+\epsilon)\, \vert{u}_x\vert^2\, dx + 
{1\over{2}} \int_{\Omega} \vert{u}_t\vert^2\, dx} +
{\lambda\over{2}} \int_0^t \int_{\Omega} \left[(\Lambda (\m)\cdot \m) u_x\right]_t dx +
 \int_0^t \int_{\Omega}  \vert{\m}_t\vert^2\, dx \\
 \\
\!\!\!\!\!\!\!\!\!\!\!\! \displaystyle{  + {1\over{4}} 
\int_{\Omega}{\left( \vert{\m}\vert^2-1\right)^2 \over \delta}\, dx + {1\over{2}} 
\int_{\Omega} {\vert{\m}_x\vert^2}\, dx }
 \displaystyle{~~\le \int_0^t \int_{\Omega}  f {u}_t\, dx \,ds + 
}\\
\\
 \displaystyle{+ {1\over{2}} 
\int_{\Omega} {\vert{\m_0}_x\vert^2}\, dx+ {1\over{2}} \int_{\Omega} 
\vert{u}_1\vert^2\, dx}\,;
\end{array}\end{equation}  
which completes the proof.
\hfill $\Box$

\bigskip
The following estimates can then be proved.
\medskip

\noindent
{\bf Lemma 2.3}
 \label{2.3}  { \it
 Let $({u}^{\varepsilon}, \m^{\varepsilon})$ denote a solution to the problem 
{\rm (\ref{eqso})--(\ref{bc-eps-del})}, then the following estimates hold{
\begin{equation}\label{ineq-lemma3}\begin{array}{cl@{\hspace{0.1ex}}c@{\hspace
{0.1ex}}l}
 \displaystyle{  \int_{\Omega} \vert{u}^{\varepsilon}_x\vert^2\, dx \le C_1} 
 \\ \\
 \displaystyle \int_{\Omega} \vert{u}^{\varepsilon}_t\vert^2\, dx \le C_2\\ \\ 
  \displaystyle \int_{\Omega} {\vert{\m}^{\varepsilon}_x\vert^2}\, dx \le C_3\\ \\
 \displaystyle
\int_{\cal Q}   \vert{\m}^{\varepsilon}_t\vert^2\, dx dt \le C_4
  \\ \\ 
 \displaystyle  \int_{\Omega}{\left( \vert{\m}^{\varepsilon}\vert^2-1\right)^2 \over \delta}\, dx \le C_5 
\end{array}\end{equation}}
where $C_k, k=1,2,3,4,5$, depend on  $ T, \m_0, f, u_1$, but do  not depend
on  $\varepsilon$ nor on $\delta$.}\rm

\medskip\noindent
{\bf Proof} {\noindent
Consider the inequality (\ref{ineq-lemma2a}) proved in Lemma 2.2, where for  simplicity,  all the supscripts
are omitted and
where the initial data are included within the constant $ C( \m_0, u_1)$,
\begin{equation}\label{ineq1b}\begin{array}{cl@{\hspace{0.1ex}}c@{\hspace
{0.1ex}}l}
 \displaystyle {1\over{2}}   \int_{\Omega} G(t+\varepsilon)\, \vert{u}_x\vert^2\, dx + 
{1\over{2}} \int_{\Omega} \vert{u}_t\vert^2\, dx + {1\over{4}} 
\int_{\Omega}{\left( \vert{\m}\vert^2-1\right)^2 \over \delta}\, dx
+ {1\over{2}}\int_{\Omega} {\vert{\m}_x\vert^2}\, dx+ 
 \\ \\
 \displaystyle{  \displaystyle \int_0^t \int_{\Omega}  \vert{\m}_t\vert^2\, dx  ~~\le -\frac{\lambda}{2}   \int_{\Omega}(\Lambda (\m)\cdot \m) u_x dx+ \int_{\Omega} \int_0^t f {u}_t\, dx \,ds + C( \m_0, u_1) }.
\end{array}\end{equation}
Now we have
\begin{equation}
 \left\vert -\frac{\lambda}{2} \int_{\Omega}(\Lambda(\m)\cdot \m) u_x dx\right\vert 
 \le \frac{\lambda}{2}  \int_{\Omega}\vert \m \vert^2 \vert u_x \vert dx.
\end{equation}
Furthermore, observe that
\begin{equation}
\int_{\Omega}\vert \m \vert^2 \vert u_x \vert dx= {\sqrt{\delta}} \int_{\Omega}{\left( \vert{\m}\vert^2-1\right) \over {\sqrt{\delta}}}  \vert u_x \vert\, dx +\int_{\Omega} \vert u_x \vert dx \le
\end{equation}
\begin{equation*}
{\sqrt{\delta}\over 2} \int_{\Omega}{\left( \vert{\m}\vert^2-1\right)^2 \over {{\delta}}}   dx +
{\sqrt{\delta}\over 2}\int_{\Omega} \vert u_x \vert^2 dx + {\sigma\over 2}\int_{\Omega} \vert u_x \vert^2 dx + {1\over {2 \sigma}} \vert \Omega\vert~,
\end{equation*}
where both $\sigma<1$ and $\delta<1$.  
Since $G(t+\varepsilon) > G(T+1), \forall t\in (0,T)$, we can choose $\sigma$ and $\delta$
so that \begin{equation}
\lambda {\sqrt{\delta}}<{1\over 2}~,  ~~~  \lambda ({\sqrt{\delta}}+ \sigma) < 
 G(T+1)
\end{equation}
and hence we obtain:
\begin{equation}\label{ineq-lemma3a}
\displaystyle{{1\over{4}}
\int_{\Omega} G(t+\varepsilon)\vert{u}_x\vert^2\, dx + {1\over{2}}
\int_{\Omega} \vert{u}_t\vert^2\, dx + \int_0^t \int_{\Omega}  \vert{\m}_t\vert^2\, dx  + {1\over{2}}
\int_{\Omega} \vert{\m}_x\vert^2\, dx + 
 }\end{equation}
{\begin{equation*}
\displaystyle{+{1\over{8}}
\int_{\Omega} {(\vert{\m}\vert^2-1)^2\over \delta}\, dx\le  \int_0^t \int_{\Omega} f u_t dx 
+ {1\over{2}} 
\int_{\Omega} {\vert{\m}_{0,x}\vert^2}\, dx +
 {1\over{2}} \int_{\Omega} \vert{u}_1\vert^2\, dx + {1\over {2 \sigma}} \vert \Omega\vert~~.}
\end{equation*}}  }
Hence, if we set {
\begin{equation}
\displaystyle E(t):= {1\over{4}}   \int_{\Omega} G(t)\, \vert{u}_x\vert^2\, dx + 
{1\over{2}} \int_{\Omega} \vert{u}_t\vert^2\, dx +{1\over{2}}\int_{\Omega} {\vert{\m}_x\vert^2}\, dx
+ {1\over{8}} \int_{\Omega}{\left( \vert{\m}\vert^2-1\right)^2 \over \delta}\, dx ~,\end{equation}  }
noting that {
    \begin{equation}
   \left\vert \int_0^t \int_{\Omega}  f {u}_t\, dx \,ds \right\vert \le {1\over{2}}\int_0^t \int_{\Omega} \vert{u}_t\vert^2\, dx \,ds +{1\over{2}}\int_0^t \int_{\Omega} \vert f \vert^2\, dx \,ds~~.
    \end{equation} }
we obtain
\begin{equation}
E(t)- \int_0^t {E}(\tau) d\tau~ \le  C(T, \m_0, f, u_1).
\end{equation}
Note, on application of Gronwall's  Lemma, it follows
that $$ E(t), \int_0^t {E}(\tau) d\tau \le  \tilde C(T, \m_0, f, u_1)~,$$
and the proof is completed since all the inequalities (\ref{ineq-lemma3}) are implied.

\hfill\eproof
}
\section{Existence result for the limit problem}

\smallskip

This Section is devoted to prove the  existence of weak solutions to the non-linear integro-differential problem 
(\ref{eq-evo}) -- (\ref{bc}). The key tools are the estimates which are independent of $\var$.
Here the limit when the  parameter $\var \to 0$   is studied. This allows us to 
establish the existence result in the generalized case of singular kernel, as far as the 
viscoelastic behaviour is concerned: this result generalizes the previous one in \cite{GVS2010}.

\bigskip
\noindent
{\bf Theorem 3.1} 
 \label{teo}  { \it  For all  $T>0$,  there exists a weak solution $(u, \m)$ to the problem  
 {\rm (\ref{eq-evo})-(\ref{ic})-(\ref{bc})},  that is a vector function $(u, \m)$  s.t. 
\begin{itemize}
\item $ u  \in L^\infty(0,T; H^1_0(\Omega))$;
\item $  u_t \in  L^\infty(0,T; L^2(\Omega))$;
\item $\m \in L^\infty(0,T; H^1(\Omega))$;
\item $ \m_t \in  L^2({\cal Q})$.
\end{itemize}
which  satisfies 
\begin{equation} \label{5-6} 
\displaystyle{ 
 \!\!\begin {array}{l} - \displaystyle \int_{\cal Q} \phi_t u^{\varepsilon}({t})  dxdt
+ 
 \int_{\cal Q} \int_0^t {G}^{\varepsilon}(t-\tau) u^{\varepsilon}_{x}(\tau)  \phi_x  d\tau  dxdt
 + 
 \int_{\cal Q}\int_0^t {\lambda\over 2}\Lambda(\m^{\varepsilon})\cdot \m^{\varepsilon} \phi_x d\tau\ dxdt
\\ \\ 
- \displaystyle  \int_{\cal Q} \left[u_1+ \int_0^t  f(\tau) d\tau\ \right] \phi dxdt
+ \displaystyle \int_{\cal Q} {\boldsymbol\psi}_t\cdot \mathbf m^{\varepsilon} dxdt
+ \int_{\cal Q} {\mathbf m_0} \cdot {\boldsymbol\psi} (\cdot, 0) dxdt
+
\\ \\
\displaystyle \int_{\cal Q} \left(\frac{|\mathbf m^{\varepsilon}|^2-1}{\delta}\right){\boldsymbol\psi} \cdot{\mathbf m^{\varepsilon}}  dxdt 
- \int_{\cal Q}  \lambda \, u^{\varepsilon}_{x} \, \Lambda (\mathbf m^{\varepsilon}) \cdot {\boldsymbol\psi} dxdt
- \int_{\cal Q}   \mathbf m^{\varepsilon}_{x} \cdot {\boldsymbol\psi}_x dxdt =0~. 
\end{array}}
\end{equation}
$\forall \phi$ smooth s.t. $ \phi(0, t)=  \phi(1, t) =0, ~ \phi(\cdot, T)=0$,  
and  $\forall {\boldsymbol\psi}\equiv(\psi_1,\psi_2)$ s.t. ~$ {\boldsymbol\psi}(x,T)=0$. 
}\rm
\medskip\noindent

{\bf Proof} \quad{\noindent
By a weak solution to
\begin{equation}\label{1}
\begin {array}{l}
\displaystyle u^{\varepsilon}_t({t})-\int_0^t {G}^{\varepsilon}(t-\tau) u^{\varepsilon}_{xx}(\tau) d\tau - u_1 
-\int_0^t \frac{\lambda}{2}(\Lambda(\m^{\varepsilon})\cdot \m^{\varepsilon})_{x} d\tau\ =\int_0^t  f(\tau) d\tau\,\\
\\  \displaystyle {u^{\varepsilon}(\cdot,0)=u_0=0,~~  u^{\varepsilon}(x, t)=0 ~ \text{\rm on} ~~ \partial\Omega\times (0,T)}
\end{array}
\end{equation}
we mean a function $u(x,t)$ such that
\begin{equation}\label{2}
\displaystyle -\int_{\cal Q} \phi_t u +\int_{\cal Q} \int_0^t {G}(t-\tau) \phi_x u_{x}(\tau) d\tau  
+ \int_{\cal Q} \int_0^t \frac{\lambda}{2}(\Lambda(\m)\cdot \m) \phi_{x} d\tau\ = 
\end{equation}
$$\displaystyle  \int_{\cal Q} \phi \left[u_1 + \int_0^t f(\tau) d\tau \right]$$
$\forall \phi$ smooth s.t. $ \phi(0, t)=  \phi(1, t) =0, ~ \phi(\cdot, T)=0$, where ${\cal Q}= \Omega\times (0,T)$ (we dropped the measure of integration $dxdt$).

By a weak solution to
\begin{equation}\label{3}
\begin {array}{l}
\displaystyle \mathbf m_t+{\mathbf m}\frac{|\mathbf m|^2-1}
{\delta}+\lambda\Lambda (\mathbf m) u_{x}-\mathbf m_{xx}=0~~~\text{\rm in} ~ {\cal Q}\\
\\  \displaystyle {{\mathbf m}(\cdot,0)={\mathbf m}_0,~~ \frac{{\partial \mathbf m}}{\partial \mbox{\boldmath
$\nu$}}=0  ~ \text{\rm on} ~~ \partial\Omega\times (0,T)}
\end{array}
\end{equation}
we mean a function ${\mathbf m}\equiv(m_1,m_2)$ such that
\begin{equation}\label{4}
 \displaystyle -\int_{\cal Q} {\boldsymbol\psi}_t \cdot \mathbf m  +  \int_{\cal Q} {\boldsymbol\psi} \cdot {\mathbf
m}\frac{|\mathbf m|^2-1}
{\delta} + \int_{\cal Q} {\boldsymbol\psi}\cdot \lambda\Lambda (\mathbf m) u_{x} + \int_{\cal Q} {\boldsymbol\psi}_x \cdot \mathbf m_{x}  \\ \\
\displaystyle \hfill -\int_{\Omega} \mathbf m_0 \cdot {\boldsymbol\psi}(x,0) dx =0 ~.
\end{equation}
$\forall {\boldsymbol\psi}$ s.t. $ {\boldsymbol\psi}(x,T)=0$ where ${\boldsymbol\psi}\equiv(\psi_1,\psi_2)$.
We start from $(u^{\varepsilon}, \m^{\varepsilon})$  satisfying 
\begin{equation} \label{5-6} 
\displaystyle{ 
\!\!\!\!\!\!\!\left \{ \!\!\begin {array}{l} - \displaystyle \int_{\cal Q} \phi_t u^{\varepsilon}({t})  
+ 
 \int_{\cal Q} \int_0^t {G}^{\varepsilon}(t-\tau) u^{\varepsilon}_{x}(\tau)  \phi_x  d\tau  
 + 
 \int_{\cal Q}\int_0^t {\lambda\over 2}\Lambda(\m^{\varepsilon})\cdot \m^{\varepsilon} \phi_x  d\tau\ =  \\
 \hskip 8cm \displaystyle
\int_{\cal Q} \left[u_1+ \int_0^t  f(\tau) d\tau\ \right] \phi  \\  
\\ 
- \displaystyle \int_{\cal Q} {\boldsymbol\psi}_t\cdot \mathbf m^{\varepsilon} 
- \int_{\cal Q} {\mathbf m_0} \cdot {\boldsymbol\psi} (\cdot, 0)
- \int_{\cal Q} \left(\frac{|\mathbf m^{\varepsilon}|^2-1}{\delta}\right){\boldsymbol\psi} \cdot{\mathbf m^{\varepsilon}}
+ \int_{\cal Q}  \lambda \, u^{\varepsilon}_{x} \, \Lambda (\mathbf m^{\varepsilon}) \cdot {\boldsymbol\psi} +\\
 \hskip 8cm \displaystyle
+ \int_{\cal Q}   \mathbf m^{\varepsilon}_{x} \cdot {\boldsymbol\psi}_x=0~. 
\end{array}\right.}
\end{equation}
\medskip
From the estimates we have 
\begin{equation}\label{7}
u^{\varepsilon}~~~~~~~\text{\rm{is bounded in}}~~~L^\infty(0,T; H^1_0(\Omega))~,
\end{equation}
\begin{equation}\label{8}
u^{\varepsilon}(\cdot,t)~~~~~~~\text{\rm{is bounded in}}~~~ H^1_0(\Omega) ~~~\text{\rm{ifor a.e. }}~t\in (0,T).
\end{equation}

Then, there exists $ u  \in L^\infty(0,T; H^1_0(\Omega))$ such that 
\begin{equation}\label{9}
u^{\varepsilon} \rightharpoonup 
~~ u~~\text{\rm{in}}~~~L^\infty(0,T; H^1_0(\Omega))~~\text{\rm{and also in}}~
L^2(0,T; H^1_0(\Omega))
\end{equation}

Moreover for a.e. $t$ there exists $v(\cdot,t)\in H^1_0(\Omega)$ such that
\begin{equation}\label{10}
u^{\varepsilon} (\cdot,t) 
\rightharpoonup ~~v(\cdot,t)~~\text{\rm{in}}~~H^1_0(\Omega), ~u^{\varepsilon} (\cdot,t) \longrightarrow v(\cdot,t)~~\text{\rm{in}}~
L^2(\Omega)~.
\end{equation}
We suppose $L^2(0,T; H^1_0(\Omega))$ equipped with the scalar product
\begin{equation}
\int_{0}^{T} \int_{\Omega} u_{x} v_{x} dx dt ~~.
\end{equation}
Let $\varphi\in{\cal D} (0,T), {\chi} \in{\cal D} ({\Omega})$, from \eqref{9} we derive
\begin{equation}\label{11}
\int_{0}^{T} \int_{\Omega} u^{\varepsilon}_{x} (\varphi {\chi})_{x} \longrightarrow
\int_{0}^{T} \int_{\Omega} u_{x} (\varphi {\chi})_{x} = 
\int_{0}^{T} \varphi \int_{\Omega} u_{x} {\chi}_{x}
\end{equation}
One has also for a.e. $t$ by \eqref{10} 
\begin{equation}\label{3.44}
\int_{\Omega} u^{\varepsilon}_{x}  {\chi}_{x} \longrightarrow
\int_{\Omega} v_{x} {\chi}_{x} ~.
\end{equation}
Moreover due to the estimates that we have $u^{\varepsilon}_{x}, \int_{\Omega} u^{\varepsilon}_{x}  {\chi}_{x}$ are uniformly bounded in $t$ and by the Lebesque theorem one derive that 
\begin{equation}\label{12}
\int_{0}^{T} \int_{\Omega} u^{\varepsilon}_{x} (\varphi {\chi})_{x} =
\int_{0}^{T} \varphi \int_{\Omega} u^{\varepsilon}_{x} {\chi}_{x} \longrightarrow
\int_{0}^{T} \varphi \int_{\Omega} v_{x} {\chi}_{x} ~.
\end{equation}
Comparing \eqref{11}, \eqref{12} we derive that for a.e. $t$
\begin{equation}\label{13}
\int_{\Omega} u^{\varepsilon}_{x} {\chi}_{x} = \int_{\Omega} v^{\varepsilon}_{x} {\chi}_{x}~~
\forall {\chi}\in {\cal D} ({\Omega})~.
\end{equation}
This implies that
\begin{equation}\label{3.47}
u (\cdot,t) = v(\cdot,t) ~~~~~\text{\rm{a.e.}} ~~ t\in (0,T).
\end{equation}
Since $ \m^{\varepsilon}$ is such that
\begin{equation}\label{14}
\m^{\varepsilon}~~~~~~~\text{\rm{is bounded in}}~~~L^\infty(0,T; H^1(\Omega))~,
\end{equation}
\begin{equation}\label{15}
\m^{\varepsilon}(\cdot,t)~~~~~~~\text{\rm{is bounded in}}~~~ H^1(\Omega)\subset C^{1/2}(\bar\Omega) ~~~\text{\rm{a.e. }}~t\in (0,T).
\end{equation}
arguing as above one derives the existence of $\m\in L^\infty(0,T; H^1(\Omega))$ s.t. up to a subsequence
\begin{equation}\label{16}
\m^{\varepsilon} \rightharpoonup
~~\m~~\text{\rm{in}}~~~L^\infty(0,T; H^1(\Omega))~,
\end{equation}
\begin{equation}\label{17}
\m^{\varepsilon} (\cdot,t) \longrightarrow \m(\cdot,t)~~\text{\rm{in}}~~C^0(\bar\Omega), ~~~\text{\rm{ a.e. }}~t.
\end{equation}
One can then pass to the limit in \eqref{5-6} to get  \eqref{2}, \eqref{4}. Let us first derive   \eqref{2}.
\begin{itemize}
\item From \eqref{9} $u^{\varepsilon} \rightharpoonup
~~ u$ in $L^2(0,T; H^1_0(\Omega))=L^2(\cal Q)$ and thus, when 
${\varepsilon}\to 0$
\begin{equation}
\int_{\Omega} {\phi}_{t} \, u^{\varepsilon}   \longrightarrow
\int_{\Omega}  {\phi}_{t} \, u ~.
\end{equation}

\item Due to our assumptions for $t$ fixed we have when ${\varepsilon}\to 0$
\begin{equation*}
 {G}^{\varepsilon}(t-\tau)    \longrightarrow ~~ {G} (t-\tau)  ~~\text{\rm{in}}~~~L^1(0,t)
\end{equation*}
and by \eqref{3.44}, \eqref{3.47}
\begin{equation*}
 \int_{\Omega} \phi_x  (x, t)\, u^{\varepsilon}_{x}(x, \tau) dx ~~ \longrightarrow ~~
  \int_{\Omega} \phi_x  (x, t)\, u_{x}(x, \tau) dx  ~~ ~~~\text{\rm{a.e. }}~\tau\in (0,t)
\end{equation*}
and these integrals are uniformly bounded independently of ${\varepsilon}$. It follows that 
\begin{equation*}
\int_0^t {G}^{\varepsilon}(t-\tau) \int_{\Omega}  \phi_x (x, t)  \, u^{\varepsilon}_{x}(x,\tau)   d\tau dx 
  \longrightarrow ~
\int_0^t {G}(t-\tau) \int_{\Omega}  \phi_x (x, t)  \, u_{x}(x,\tau)   d\tau dx\,.
\end{equation*}
Since these integrals are  uniformly bounded independently  it follows that 
\begin{equation*}
\int_0^T\!\!\!\!\!\! \int_0^t {G}^{\varepsilon}(t-\tau)\!\! \int_{\Omega}  \phi_x (x, t)  \, u^{\varepsilon}_{x}(x,\tau)   d\tau dx dt
  \rightarrow ~
\int_0^T\!\!\! \!\!\!\int_0^t {G}(t-\tau) \!\! \int_{\Omega}  \phi_x (x, t)  \, u_{x}(x,\tau)   d\tau dx  dt
\end{equation*}
i.e.
\begin{equation*}
\int_{\cal Q}  \int_0^t {G}^{\varepsilon}(t-\tau)   \, u^{\varepsilon}_{x} \,  \phi_x
  \longrightarrow ~
\int_{\cal Q}  \int_0^t {G}(t-\tau)\, u_{x}   \,  \phi_x~.
\end{equation*}
\item Since $\mathbf m^{\varepsilon} (\cdot,\tau) \to \mathbf m (\cdot,\tau)$ in $C^0(\bar\Omega)$ a.e. $\tau$ one has 
\begin{equation*}
 \int_{\cal Q} \int_0^t  {\lambda\over 2}  \Lambda (\mathbf m^{\varepsilon}) \cdot \mathbf m^{\varepsilon}\,\phi_{x} d\tau  \longrightarrow 
  \int_{\cal Q} \int_0^t  {\lambda\over 2}  \Lambda (\mathbf m) \cdot \mathbf m \,\,\phi_{x} d\tau
\end{equation*}
which completes the existence of solution to \eqref{2}.
To pass to the limit in the second equation in \eqref{5-6} due to \eqref{16}-\eqref{17} only perhaps the fourth integral is not clear. But from 
\eqref{17}, \eqref{10}
\begin{equation*}
\Lambda (\mathbf m^{\varepsilon}) \cdot {\boldsymbol\psi} (\cdot,t) \longrightarrow 
\Lambda (\mathbf m) \cdot {\boldsymbol\psi} (\cdot,t)~~\text{\rm{in}}~~L^2(\Omega), ~~~~~~~\text{\rm{ a.e. }}~t.
\end{equation*}
\begin{equation*}
u^{\varepsilon}_x (\cdot,t) \rightharpoonup
~~ u_x(\cdot,t)~~~~\text{\rm{in}}~~~ L^2(\Omega), ~~~~~
~~~~~~~\text{\rm{a. e.}}~~~ t.
\end{equation*}
This implies that 
\begin{equation*}
\int_{\Omega} \lambda \,  u^{\varepsilon}_x  \Lambda (\mathbf m^{\varepsilon}) \cdot {\boldsymbol\psi} (\cdot,t) \longrightarrow 
\int_{\Omega} \lambda \, u_x \Lambda (\mathbf m) \cdot {\boldsymbol\psi} (\cdot,t)~~ ~~~~~~~\text{\rm{ a.e. }}~t.
\end{equation*}
Since all these integrals are uniformly bounded one deduces that
\begin{equation*}
\int_{\cal Q} \lambda \,  u^{\varepsilon}_x \, \Lambda (\mathbf m^{\varepsilon}) \cdot {\boldsymbol\psi} \longrightarrow 
\int_{\cal Q} \lambda \, u_x \, \Lambda (\mathbf m) \cdot {\boldsymbol\psi} 
\end{equation*}
which completes the proof of the existence of a solution to \eqref{4}. 
\end{itemize}

\hfill\eproof}

\begin{center} {\bf Acknowledgments} \end{center}
S.Carillo wishes to acknowledge the partial financial  support of GNFM-INDAM, INFN and SAPIENZA 
Universit\`a di Roma. The research of M. Chipot leading to these results has received funding from Lithuanian-Swiss cooperation programme to reduce economic and social disparities within the enlarged European Union under project agreement 
No CH-3-SMM-01/0.
M. Chipot  thanks also the Italian GNFM-INdAM and SBAI Dept. for the kind hospitality in Rome.

\end{document}